\documentclass[a4paper,notitlepage, 8pt,reqno]{article}

\usepackage[cp1251]{inputenc}
 \usepackage[T2A]{fontenc}
 \usepackage[russian,english]{babel}
\usepackage[all,arc,poly,2cell,curve,arrow,tips]{xy}
\usepackage{enumerate, eucal, amsthm,amsmath, amssymb}
\usepackage{graphicx}
\usepackage{mathtext}

\usepackage{tikz}

\allowdisplaybreaks[4]

\theoremstyle{definition}

\theoremstyle{plain}
\newtheorem{thm}{Theorem}
\newtheorem*{thm*}{Theorem}
\newtheorem{prop}{Proposition}
\newtheorem*{prop*}{Proposition}
\newtheorem{cor}{Corrolary}
\newtheorem{lem}{Lemma}

\title{ A Gelfand-Tsetlin type base for the algebra $\mathfrak{sp}_4$ and hypergeometric functions}

\author{D.V. Artamonov\footnote{Lomonosov Moscow State University, artamonov.dmitri@gmail.com}}

\begin{document}
\maketitle


In the paper a realization of representation of a Lie algebra  $\mathfrak{sp}_4$  in the space of function on the Lie group  $Sp_4$ is considered.
We  find a function corresponding to a Gelfand-Tsetlin type vector  for  $\mathfrak{sp}_4$ constructed by D.P. Zhelobenko.  This function is expressed though a  $A$-hypergeometric function.  After developing some new technique we derive analytically formulas for the action of generators of algebra in this base (the were not known before).  These formula turn out to be much more difficult than the formulas for the action of generators in the Gelfand-Tsetlin type base constructed by Molev.



\section{Introduction}

In the paper Lie groups and algebras over   $\mathbb{C}$ are considered.Consider a Lie group  $Sp_4$  and functions on it.  The group acts on the space of function on itself by right  shifts and the corresponding Lie algebra  acts by infinitesimal right shifts.
Finite dimensional irreducible representation of  the Lie algebra $\mathfrak{sp}_4$ can be realized as subrepresentations in this functional representation  \cite{zh2}.

Also other construction of finite-dimensional irreducible representations of  $\mathfrak{sp}_4$ are known.  Thus there exists a construction of Gelfand and Tsetlin. In the book \cite{zh2}  Zhelobenko gave a construction of base vectors of Gelfand-Tsetlin type  for a representation of $\mathfrak{sp}_{2n}$. This construction is based on equivalence of restriction problems   $\mathfrak{gl}_{n+1}\downarrow\mathfrak{gl}_{n-1}$  and  $\mathfrak{sp}_{2n}\downarrow\mathfrak{sp}_{2n-2}$.  Here by a problem of restriction   $\mathfrak{g}\downarrow\mathfrak{h}$,     where  $\mathfrak{g}\supset\mathfrak{h}$ are a  Lie algebra and a subalgebra in it, we mean a problem of description of the space of     $\mathfrak{h}$-highest vectors with a fixed  $\mathfrak{h}$-weight in a given irreducible representation of  $\mathfrak{g}$.  An equivalence of these problems is just a linear isomorphism between the corresponding linear spaces. Later Molev constructed  (using other ideas  \cite{M})  a Gelfand-Tsetlin type base for a representation of   $\mathfrak{sp}_{2n}$.
Moreover, Molev derived formulas for the action of generators of  the algebra in this base.  Such formulas for the Zhelobenko's base were not known.

Let us return to the Zhelobenko's base. There appeares a question. What is a relation between the Zhelobenko's construction and the functional realization.? More precise, {\it which function on a group there corresponds to a Gelfand-Tsetlin type vector?}

This question is natural since the Zhelobenkos's construction uses the functional realization. There exist two cases when an answer to a similar question is known. These are the cases of the Gelfand-Tsetlin bases for the algebras $\mathfrak{gl}_2$ and  $\mathfrak{gl}_3$.  The case of  $\mathfrak{gl}_2$ is trivial and the answer in the case  $\mathfrak{gl}_3$ is both nontrivial and very beautiful. The function corresponding to a Gelfand-Tsetlin vector can be expressed through the Gauss hypergeometric function  \cite{bb}.  From  a modern viewpoint it is more natural to express the function thought an  $A$-hypergeometric function \cite{a1}.

The explicit formulas for functions corresponding to Gelfand-Tsetlin  vectors are used in  \cite{bb}  to obtain formulas for the action of generators of the algebra in this base.  There derivation uses complex analysis and the theory of special functions.  In 
 \cite{A1}    they are used for derivation of explicit formulas for Clebsh-Gordon coefficients for the algebra  $\mathfrak{gl}_3$.   They can be used for obtaining explicit constructions of infinite-dimensional representations  of  $\mathfrak{gl}_3$  \cite{tmph}.  There exist their generalizations to the case of quantized algebras   \cite{t15}, \cite{t16}.

 In the paper we consider the question in the case of the algebra  $\mathfrak{sp}_4$.   From one hand  in this case we manage to derive a fomula for a function corresponding to a Gelfand-Tsetlin-Zhelobenko vector  and the answer is both nontrivial and quite simple. From the other hand the algebra  $\mathfrak{sp}_4$ is of special interest among  symplectic algebras from both mathimatical point of view  (as the simplest example of the series   $C$)  \cite{ap1},  \cite{ap2} and also form the point of view of  physics \cite{ap3}, \cite{ap4}.

The main result of the paper is the following. Using a relation between the restriction problems   $\mathfrak{sp}_4\downarrow\mathfrak{sp}_2$  and  $\mathfrak{gl}_3\downarrow\mathfrak{gl}_1$ we derive that a function corresponding to a Gelfand-Tsetlin-Zhelobenko base vector  can be expressed though a  $A$-hypergeometric function. But in contrast to the case $\mathfrak{gl}_3$ this function cannot be reduced  to a function of one argument.

The fact that a function corresponding to a Gelfand-Tsetlin-Zhelobenko base vector can be expressed  through a  $A$-hypergeometric function is a notable fact that does not take place for the algebras  $\mathfrak{gl}_n$, $\mathfrak{sp}_{2n}$ for bigger $n$.

Then we use these formulas to derive formulas for the action of generators  of  $\mathfrak{sp}_4$  in the  Gelfand-Tsetlin-Zhelobenko base.
This result was not obtained by Zhelobenko.  It is interesting that these formulas are much more difficult than the formulas for the action of generators in the Molev's base.

Using the functional approach  and the theory of   $A$-hypergeometric functions we manage to obtain quite simple formulas for the coefficients in the formulas for the action.

Let us note that to derive the formulas for the action of generators we develop some new technique  dealing with  $A$-hypergeometric functions and introduce some new classes of hypergeometric type functions..

\section{The basic notions}

\subsection{The algebra  $\mathfrak{sp}_{4}$}
\label{algl}

The Lie algebra  $\mathfrak{sp}_{4}$ is considered as a subalgebra in the Lie algebra of all  $4\times 4$ matrices. We index rows and columns by labels   $i,j=-2,-1,1,2$.
The algebra  $\mathfrak{sp}_{4} $  is spanned by matrices
\begin{equation}\label{fd}F_{i,j}=E_{i,j}-sign(i)sign(j)E_{-j,-i},\end{equation}
where  $i,j=-2,-1,1,2$. The subalgebra   $\mathfrak{sp}_{2} $ is generated by
 $<F_{i,j}>_{i,j=-2,2}$.

We take a Lie algebra  $\mathfrak{gl}_{3}$ of all  $3\times 3$ matrices. We index rows and columns by labels    $i,j=-2,-1,1$. This algebra is spanned by matrices  $E_{i,j}$,   $i,j=-2,-1,1$.  Chose a subalgebra   $\mathfrak{gl}_{1}=<E_{i,j}>_{i,j=-2}$.

\subsection{Functions on a group}
\subsubsection{A functional realization of a representation}
\label{fugu}
We need a realization of representation on the space of functions on a group $G=Sp_{4}$,  $GL_{3}$ (see  \cite{zh2}). Onto a function
$f(g)$, $g\in G$ an element $X\in G_{}$ acts by a right shift

\begin{equation}
\label{xf}
(Xf)(g)=f(gX).
\end{equation}

Fix a highest weight   $[m]=[m_{-2},m_{-1}]$  in the case  $Sp_{4}$ and  $[m]=[m_{-2},m_{-1},m_{1}=0]$ in the case  $GL_{3}$.


\subsubsection{Determinants}
Let us give an example of a function on a group. Let  $a_{i}^{j}$ be a function of a matrix element on  $GL_{3}$ or  $Sp_4$. Here  $j$ is a row label and  $i$is a column label. Put

\begin{equation}
\label{dete}
a_{i}:=a_i^{-2},\,\,\,\, a_{i_1,i_2}:=det(a_i^j)_{i=i_1,i_2}^{j=-2,-1}.
\end{equation}


Using formulas   \eqref{xf}, \eqref{dete}  we obtain that onto  $a_{i}$  and   $a_{i_1,i_2}$ the operator  $E_{i,j}$
acts through a transformation of column labels by the ruler

\begin{equation}
\label{edet1}
E_{i,j}a_{i_1,}=a_{\{i_1\}\mid_{j\mapsto i}},\,\,\,\,\,E_{i,j}a_{i_1,i_2}=a_{\{i_1,i_2\}\mid_{j\mapsto i}},
\end{equation}

where  $.\mid_{j\mapsto i}$ is an operation of substitution of  $j$ instead of
$i$,  if  $j$ is one of the indices of   $a$. One  obtains $0$ otherwise.  An operator  $F_{i,j}$
can be expressed by formula  
\eqref{fd}.


\subsubsection{Formulas for the highest vectors}

Let us give an explicit formula for a highest vector of weight $[m_{-2},m_{-1},0]$ for $\mathfrak{gl}_3$ and of weight  $[m_{-2},m_{-1}]$  for $\mathfrak{sp}_4$. The formula is the following one:
\begin{equation}\label{vno}v_0=a_{-2}^{m_{-2}-m_{-1}}
a_{-2,-1}^{m_{-1}}
\end{equation}

 For the algebra   $\mathfrak{gl}_3$ all  possible $\mathfrak{gl}_1$-highest vectors form a span of vectors   $E_{-1,-2}^aE_{1,-1}^bv_0$.  For the algebra  $\mathfrak{sp}_4$ all  possible  $\mathfrak{sp}_2$-highest vectors form a span of vectors $F_{-1,-2}^aF_{1,-1}^bv_0$. Using the formula  \eqref{edet1}
we obtain the Lemma.

\begin{lem}
\label{le5}
For the algebra   $\mathfrak{gl}_3$  the  $\mathfrak{gl}_1$-highest vectors  can be expressed as polynomials in  $a_{-2}, a_{\pm 1},a_{-2,\pm1},a_{-1,1}$.

For the algebra  $\mathfrak{sp}_4$ the     $\mathfrak{sp}_2$-highest vectors  can be expressed as polynomials in  $a_{-2}, a_{\pm 1},a_{-2,\pm1},a_{-1,1}$.

\end{lem}

\subsection{The Gelfand-Tsetlin base}

One can find the   details in \cite{zh2}.

A Gelfand-Tsetlin base for a  chain of algebras  $\mathfrak{gl}_3 \supset\mathfrak{gl}_2\supset\mathfrak{gl}_1$  is a base that is obtained as follows. Let $V_{\mu_3}$  be an irreducible finite-dimensional representation of   $\mathfrak{gl}_3 $  with the highest weight  $\mu_3$.  Consider it as a representation of  $\mathfrak{gl}_2 $,   then it splits into a sum of  $\mathfrak{gl}_2 $-irreducible representations.  It is known that (see \cite{zh2}) an irreducible representation of  $\mathfrak{gl}_2 $ with a given highest weight  $\mu_2$ can occur in this sum only with a multiplicity one. Denote it as  $V_{\mu_3,\mu_2}$.  thus we obtain

 $$
V_{\mu_3}=\oplus_{\mu_{2}}  V_{\mu_3,\mu_2},
$$

Now consider $ V_{\mu_3,\mu_2}$ as a representation of $\mathfrak{gl}_1$   and split it into   $\mathfrak{gl}_1$-irreducible representations. An irreducible representation of  $\mathfrak{gl}_1$ with a given weight  $\mu_1$ can occur in this sum only with a multiplicity one. Denote it as $V_{\mu_3,\mu_2,\mu_1}$.
One obtains

 $$
V_{\mu_3}=\oplus_{\mu_{2}} \oplus_{\mu_1} V_{\mu_3,\mu_2,\mu_1}.
$$

since irreducible representations of  $\mathfrak{gl}_1$ are  one-dimensional then taking   a vector in each $ V_{\mu_3,\mu_2,\mu_1}$  we obtain a base $V_{\mu_3}$, which is called the Gelfand-Tsetlin base. The base vectors are encoded by a collection of highext vectors  $\mu_3,\mu_2,\mu_1$.  If one writes them one under another one gets an integer  diagram

\begin{equation}
	\label{dgc3}
\begin{pmatrix}
m_{-2} && m_{-1} &&0\\ &k_{-2}&& k_{-1}\\&&h_{-2}
\end{pmatrix},
\end{equation}

in which  the betweeness conditions hold. That is if one element is written between two elements of a higher row then if lies between them. The inverse is true: every integer diagram for which the betweeness condition holds appears as a Gelfand-Tsetlin diagram for a representation with the highest weight  $[m_{-2},m_{-1},0]$.

For a chain $\mathfrak{sp}_4\supset\mathfrak{sp}_2$ a Gelfand-Tsetlin type base is a construction of the following type. Let  $V_{\mu_4}$ be a finite dimensional representation of    $\mathfrak{sp}_4 $  with the highest weight  $\mu_4$.   Consider it as a representation of   $\mathfrak{sp}_2 $.    It splits into a sum of  $\mathfrak{sp}_2 $-irreducible representations. But an irreducible  representation  $V_{\mu_4, \mu_2} $ of   $\mathfrak{sp}_2 $ with a fixed highest  weight   $\mu_2$ can occur with some multiplicity. Thus we write

$$
V_{\mu_4}=\oplus_{\mu_{2}}  Mult_{\mu_2} \otimes V_{\mu_4, \mu_2},
$$

where  $Mult_{\mu_2}$ is a vector space of   $\mathfrak{sp}_2$-highest vectors with a fixed weight.  Let   $\mu'_4$  be a label indexing some base in     $Mult_{\mu_2}$.
Since   $\mathfrak{sp}_2\simeq\mathfrak{sl}_2$,  then in a representation $V_{\mu_2}$  indexed by  $\mu_4,\mu'_4,\mu_2,\mu'_2$.  Such a base is called a Gelfand-Tsetlin type base.  Since  a base  $Mult_{\mu_2}$ was not fixed, then the  Gelfand-Tsetlin type base is mot unique. 

If one writes   $\mu_4,\mu'_4,\mu_2,\mu'_2$ one under another then one obtains a Gelfand-Tsetlin diagram. It's structure will be explained below.

\subsection{  $A$-hypergeometric function}

\subsubsection{A $\Gamma$-series}  A detailed information can be found in  \cite{GG}.

Let  $B\subset \mathbb{Z}^N$ be a lattice and let $\gamma\in \mathbb{Z}^N$ be a fixed vector.  Define a  {\it hypergeometric
	$\Gamma$-series  } in variables $z_1,...,z_N$ by the formula

\begin{equation}
\mathcal{F}_{\gamma}(z)=\sum_{b\in
	B}\frac{z^{b+\gamma}}{\Gamma(b+\gamma+1)},
\end{equation}
where  $z=(z_1,...,z_N)$, and we use the notations

$$
z^{b+\gamma}:=\prod_{i=1}^N
z_i^{b_i+\gamma_i},\,\,\,\Gamma(b+\gamma+1):=\prod_{i=1}^N\Gamma(b_i+\gamma_i+1).
$$

Note that the set of powers   of summands in a $\Gamma$-series form a shifted lattice $\gamma+B$.

Note that if for a summand in a  $\Gamma$-series one of the numbers $b_i+\gamma_i$ negative then this summand vanishes.

Let us formulate some properties of a  $\Gamma$-series:

\begin{enumerate}

	\item A vector $\gamma$ can be changes to $\gamma+b$, $b\in B$,  the series does not change
	
	\item A $\Gamma$-series can be differentiated by the ruler:   \begin{equation}
	\label{dgr}\frac{\partial }{\partial
		z_i}\mathcal{F}_{\gamma}(z)=\mathcal{F}_{\gamma-e_i}(z),\end{equation} where  $e_i=(0,...,1,...,0)$,  where
	  $1$ occurs at the place  $i$.

	\item
	Let  $F_{2,1}(a_1,a_2,b_1;z)=\sum_{n\in\mathbb{Z}^{\geq
			0}}\frac{(a_1)_n(a_2)_n}{(b_1)_n}z^n$, where
	$(a)_n=\frac{\Gamma(a+n)}{\Gamma(a)}$, be a Gauss'
	hypergeometric series. Then if $\gamma=(-a_1,-a_2,b_1-1,0)$,  and
	$B=\mathbb{Z}<(-1,-1,1,1)>$, then
	
	\begin{align*}
	&\mathcal{F}_{\gamma}(z_1,z_2,z_3,z_4)=cz_1^{-a_1}z_2^{-a_2}z_3^{b_1-1}F_{2,1}(a_1,a_2,b_1;\frac{z_3z_4}{z_1z_2})\\
	&c=\frac{1}{\Gamma(1-a_1)\Gamma(1-a_2)\Gamma(b_1)}
	\end{align*}

\end{enumerate}

A sum of a  $\Gamma$-series  (if it converges) is called a   $A$-hypergeometric function.

A $A$-hypergeometric function satisfies a system of PDE which is called the  Gelfand-Kapranov-Zelevinsky system (GKZ shortly). It consists of equations of two types.

{\bf 1.}  Let $a=(a_1,...,a_N)$  be a vector orthogonal to  $B$,  then

\begin{equation}
\label{e1}
a_1z_1\frac{\partial}{\partial z_1}\mathcal{F}_{\gamma}+...+a_Nz_N\frac{\partial}{\partial z_N}\mathcal{F}_{\gamma}=(a_1\gamma_1+...+a_N\gamma_N)\mathcal{F}_{\gamma},
\end{equation}
It is sufficient to consider only base vectors in the orthogonal complement to  $B$.

{\bf 2.} Let  $b\in B$ and $b=b_+-b_-$,  where coordinates of  $b_+$, $b_-$ are non-negative.  Take non-zero elements in these vectors
$b_+=(...b_{i_1},....,b_{i_k}...)$,  $b_-=(...b_{j_1},....,b_{j_l}...)$.  Then

\begin{equation}
\label{e2} (\frac{\partial }{\partial
z_{i_1}})^{b_{i_1}}...(\frac{\partial}{\partial z_{i_k}})^{b_{i_k}}
\mathcal{F}_{\gamma}=(\frac{\partial }{\partial
z_{j_1}})^{b_{j_1}}...(\frac{\partial }{\partial z_{j_l}})^{b_{j_l}} \mathcal{F}_{\gamma}
\end{equation}

\subsection{The case $\mathfrak{gl}_3$}
\label{ra2}

Let us give a formula for a function corresponding to a diagram   \eqref{dgc3} for  $\mathfrak{gl}_3$.
Such a formula is given in the next  Theorem proved in \cite{bb}.

\begin{thm}\label{vec3} 
	
  Take the determinants in the following order  $$a=(a_{-2},a_{-1},a_{1},a_{-2,-1},a_{-2,1},a_{-1,1}),$$
	
	take a lattice
	
	$$
B=\mathbb{Z}<	(1,-1,0,0,-1,1)>.
	$$

	 $\gamma=(h_{-2}-m_{-1},k_{-2}-h_{-2},    m_{-1}-k_{-1}        ,        k_{-2}    ,m_{-1}-k_{-1},0)$. Then to a diagram there corresponds a function  $\mathcal{F}_{\gamma}(a)$
\end{thm}

A more explicit formula for the function  $\mathcal{F}_{\gamma}(a)$  is given in \eqref{ac} below.

The obtained  $\Gamma$-series be expressed through a Gauss' hypergeometric series.  In this form this Theorem was obtained in \cite{bb}.


Note that the lattice   $B$  can be defined by equations onto powers of determinants:

$$
\begin{cases}
\text{ the sum of powers of determinants that contain indices  $-2$, or $-1$,or $1$}=m_{-2},\\
\text{  the sum of powers of determinants that contain indices   $-2$ and $-1$, $-2$ and $1$, $-1$ and  $1$}=m_{-1}\\
\text{  the sum of powers of determinants that contain indices   $-2$ or  $-1$ }=k_{-2}\\
\text{  the sum of powers of determinants that contain indices   $-2$ and $-1$ }=k_{-1}\\
\text{  the sum of powers of determinants that contain indices   $-2$}=h_{-2}
\end{cases}
$$





Note that to a Gelfand-Tsetlin diagram there corresponds a shifted lattice. To a shifted lattice there corresponds a shift  $\gamma$, defined    $modB$.

\section{A function corresponding to a Gelfand-Tsetlin vector in the case  $\mathfrak{sp}_4$}
\label{ac2}
Let us be given a representation $\mathfrak{sp}_4$ with the highest weight
$[m_{-2},m_{-1}]$. consider a restriction problem
$\mathfrak{sp}_4\downarrow\mathfrak{sp}_2$. In   \cite{zh2} it is shown that  the  restriction problems
$\mathfrak{gl}_3\downarrow\mathfrak{gl}_1$ and
$\mathfrak{sp}_{4}\downarrow\mathfrak{sp}_{2}$ ar equivalent. Explicitly to an expression in   $a_i^j$, considered as function on $GL_3$ there coresponds the same expression in   $a_i^j$ considered as function on  $Sp_4$.
Thus,
$\mathfrak{sp}_2$-highest vectors are encoded by  integer diagrams  \eqref{dgc3}. 

To a diagram there corresponds a function $\mathcal{F}_{\gamma}(a)$,  which can be written explicitly as follows (we use a notation $(x_1,...,x_n)!:=x_1!\cdot...\cdot x_n!$):
\begin{equation}
	\label{ac} \frac{a_{1}^{m_{-2}-k_{-2}}a_{-2,-1}^{k_{-2}}}{(m_{-2}-k_{-2})!k_{-2}!}\sum
	\frac{	a_{-1}^{p_{-1
		}}a_{-1,1}^{p_{-1,1}}a_{-2}^{p_{-2}}a_{-2,1}^{p_{-2,1}}}{p_{-1}!p_{-1,1}!p_{-2}!p_{-2,1}!}
,
\end{equation}
where a summation is taken over all positive  integers $p_{-1},p_{-1,1},p_{1},p_{-2,1}$, such that

\begin{equation}
	p_{-1}+p_{-2}=k_{-2}-m_{-1},\,\,\,p_{-1,1}+p_{-2,1}=m_{-1}-k_{-1},\,\,\, p_{-1}+p_{-1,1}=k_{-2}-h_{-2}.
\end{equation}

This sum is actually finite.

Note that  \eqref{dgc3} is a part  $\mu_4,\mu'_4,\mu_2$ of a diagram for  $\mathfrak{sp}_4$.

Thus we obtained a formula for a   $\mathfrak{sp}_2$-highest vector. It's highest weight is     $[h_{-2}]$.  
Let us use that $\mathfrak{sp}_2\simeq \mathfrak{sl}_2$.  A standard base in a representation of   $\mathfrak{sl}_2$  of the highest weight $[h_{-2}]$ is encoded by diagrams
\begin{equation}
\label{dsp2}
 \begin{pmatrix}
 h_{-2} && 0\\& h_{-1}
 \end{pmatrix},
\end{equation}
 
and a vector corresponding to this diagram can be obtained form the highest vector by applying   of an operator $\frac{F_{2,-2}^{h_{-2}-h_{-1}}}{(h_{-2}-h_{-1})!}$.
If one applies it to  \eqref{ac} then one gets

\begin{equation}
\label{ac1} \frac{a_{1}^{m_{-2}-k_{-2}}}{(m_{-2}-k_{-2})!}\sum
\frac{a_{-2,-1}^{k'_{-2}}a_{2,-1}^{k''_{-2}}a_{-1}^{p_{-1
	}}a_{-1,1}^{p_{-1,1}}a_{-2}^{p'_{-2}}a_{2}^{p''_{-2}}a_{-2,1}^{p'_{-2,1}}   a_{2,1}^{p''_{-2,1}} }{k'_{-2}k''_{-2}p_{-1}!p_{-1,1}!p'_{2}!p''_{2}!p'_{-2,1}!p''_{-2,1}!}
,
\end{equation}

where a summation is taken over all positive  integers $p_{-1},p_{-1,1},p_{1},p_{-2,1}$,  such that

\begin{equation}
\begin{cases}
k'_{-2}+k''_{-2}=k_{-2},,\,\,\, p'_{-2}+p''_{-2}=p_{-2},\,\,\, p'_{-2,1}+p''_{-2,1}=p_{-2,1},\\
k''_{-2}+p''_{-2}+p''_{-2,1}=h_{-2}-h_{-1}\\
p_{-1}+p_{-2}=k_{-2}-m_{-1},\,\,\,p_{-1,1}+p_{-2,1}=m_{-1}-k_{-1},\,\,\, p_{-1}+p_{-1,1}=k_{-2}-h_{-2}.
\end{cases}
\end{equation}

These equations can be described as follows

\begin{equation}
\label{ur}
\begin{cases}
\text{ the sum of powers of determinants that contain indices $\pm 2$, or $-1$, or $1$}=m_{-2},\\
\text{ the sum of powers of determinants that contain indices $\pm 2$  and  $-1$, $\pm 2$ and  $1$, $-1$ and $1$}=m_{-1}\\
\text{ the sum of powers of determinants that contain indices $\pm 2$ or $-1$ }=k_{-2}\\
\text{ the sum of powers of determinants that contain indices $\pm 2$ and $-1$ }=k_{-1}\\
\text{ the sum of powers of determinants that contain indices $\pm 2$ }=h_{-2}\\
\text{ the sum of powers of determinants that contain indices  $-2$} =h_{-1}
\end{cases}
\end{equation}

As before to a diagram there corresponds a shifted lattice. To a shifted lattice there corresponds a shift  $\gamma$ defined   $modB$.

There exist a graphical encoding of these equations, it is given in a picture \eqref{dgu} and the text after it.

Everywhere below we denote as  $B$ we denote a lattice in the space  with coordinates whose coordinates are labeled by the determinants

 $$(a_{-2},a_{-1},a_{1},a_{2},a_{-2,-1},a_{-2,1},a_{-2,2},a_{-1,1},a_{-1,2},a_{1,2})$$


Let us write a base in   $B$.   It consists of three vectors  $v_1,v_2,v_3$,    written  as rows follows

\begin{align}
\begin{pmatrix}
a_{-2}&a_{-1}&a_{1}&a_{2}&a_{-2,-1}&a_{-2,1}&a_{-2,2}&a_{-1,1}&a_{-1,2}&a_{1,2}\\
1&-1&0&0&0&-1&0&1&0&0\\
-1&0&0&1&1&0&0&0&-1&0\\
-1&0&0&1&0&1&0&0&0&1
\end{pmatrix}
\end{align}

By concatenating  \eqref{dgc3} and  \eqref{dsp2}we obtain that a Gelfand-Tsetlin diagram for  $\mathfrak{sp}_4$ is encoded by

\begin{equation}
\label{dsp4}
\begin{pmatrix}
m_{-2} && m_{-1} &&0 &&\\&k_{-2}&& k_{-1} &&\\&&h_{-2}&&0\\&&&h_{-1}
\end{pmatrix}
\end{equation}

Thus we have proved a Theorem

\begin{thm}
To a Gelfand-Tsetlin diagram for  $\mathfrak{sp}_4$ of type  \eqref{dsp4} there corresponds a  $\Gamma$-series in determinants  $a_{\pm 2}, a_{\pm 1}, a_{\pm 2,\pm 1}, a_{-1,1}$ defined by a shifted lattice   \eqref{ur}.
\end{thm}


\section{The action of generators of the algebra}

It is enough to describe an action of generators  $F_{-2,-2}$ and $F_{-1,-1}$,  $F_{-2,2}$ and  $F_{2,-2}$,  $F_{-2,1}$ and   $F_{-1,-2}$,  .

\subsection{Operators  $F_{-2,-2}$ and  $F_{-1,-1}$}

When  $F_{-2,-2}=E_{-2,-2}-E_{2,2}$  acts onto a product of determinants, this product is multiplied onto  a difference of the number of occurrences of indices    $-2$ and  $2$  in these determinants. Due to   \eqref{ur}  this number is the same for all summands in our   $\Gamma$-series and it equals to  $h_{-2}-h_{-1}$.  Thus the vector  \eqref{dsp4} is an eigenvector for   $F_{-2,-2}$ with an eigenvalue   $h_{-2}-h_{-1}$.

When   $F_{-1,-1}=E_{-1,-1}-E_{1,1}$  acts onto a product of determinants, this product is multiplied onto  a difference of the number of occurrences of indices    $-1$  and  $1$   in these determinants. Due to  \eqref{ur}  this number is the same for all summands in our    $\Gamma$-series and it equals to  $2(k_{-2}+k_{-1})-(m_{-2}+m_{-1})-h_{-2}$. Thus the vector \eqref{dsp4} s an eigenvector for  $F_{-1,-1}$ with an eigenvalue  $2(k_{-2}+k_{-1})-(m_{-2}+m_{-1})-h_{-2}$.

\subsection{Operators  $F_{-2,2}$ and  $F_{2,-2}$}

By definition   \eqref{dsp4} is obtained form a   $\mathfrak{sp}_2$-highest vector by applying  of the operator $\frac{F_{2,-2}^{h_{-2}-h_{-1}}}{(h_{-2}-h_{-1})!}$.

Hence after applying  $F_{2,-2}$ we obtain a diagram which is obtained from an original diagram by the trasformation   $h_{-1}\mapsto h_{-1}-1$,  taken with a coefficient    $h_{-2}-h_{-1}+1$.

And applying the operator  $F_{-2,2}$  we obtain a diagram which is obtained from an original diagram by the  trasformation   $h_{-1}\mapsto h_{-1}+1$, taken with a coefficient    $h_{-1}+1$.

\subsection{Operators  $F_{-2,1}$ и  $F_{1,-2}$. Differential operators}

The considered operators can be writen as the following differential operators

\begin{align*}
& F_{-2,1}=a_{-2}\frac{\partial}{\partial a_{1}}+a_{-2,-1}\frac{\partial}{\partial a_{1,-1}}+a_{-2,2}\frac{\partial}{\partial a_{1,2}}+
a_{-1}\frac{\partial}{\partial a_{2}}+a_{-2,-1}\frac{\partial}{\partial a_{-2,2}}+a_{-1,1}\frac{\partial}{\partial a_{2,1}}=\\
&=a_{-2}\frac{\partial}{\partial a_{1}}+
a_{-1}\frac{\partial}{\partial a_{2}}-a_{-2,-1}\frac{\partial}{\partial a_{-1,1}}-2 a_{-1,1}\frac{\partial}{\partial a_{1,2} },\\
&F_{1,-2}=a_{1}\frac{\partial}{\partial a_{-2}}+a_{1,-1}\frac{\partial}{\partial a_{-2,-1}}+a_{1,2}\frac{\partial}{\partial a_{-2,2}}+
a_{2}\frac{\partial}{\partial a_{-1}}+a_{-2,2}\frac{\partial}{\partial a_{-2,-1}}+a_{2,1}\frac{\partial}{\partial a_{-1,1}}=\\
&=a_{1}\frac{\partial}{\partial a_{-2}}+
a_{2}\frac{\partial}{\partial a_{-1}}-a_{1,2}\frac{\partial}{\partial a_{-1,1}}-2a_{-1,1}\frac{\partial}{\partial a_{-2,-1}}.
\end{align*}

Note that  $a_{-1,1}=-a_{-2,2}$.
Due to  \eqref{dgr}, one has

\begin{align}
\begin{split}
\label{ffd}
&F_{-2,1}\mathcal{F}_{\gamma}=a_{-2}\mathcal{F}_{\gamma-e_{1}}+a_{-1}\mathcal{F}_{\gamma-e_{2}}-a_{-2,-1}\mathcal{F}_{\gamma-e_{-1,1}}-2a_{-1,1}\mathcal{F}_{\gamma-e_{1,2}}
,\\
&F_{1,-2}\mathcal{F}_{\gamma}=a_{1}\mathcal{F}_{\gamma-e_{-2}}+a_{2}\mathcal{F}_{\gamma-e_{-1}}-a_{1,2}\mathcal{F}_{\gamma-e_{-1,1}}-2a_{-1,1}\mathcal{F}_{\gamma-e_{-2,-1}}
\end{split}
\end{align}

To  obtain  explicit formulas for the action of  $F_{-2,1}$ and  $F_{1,-2}$ we need formulas for a product of an $A$-hypergeometric function and a variable that hold modulo the Plucker relations.

\subsection{ Functions  $\mathcal{F}_{\gamma}^s$ and   $F_{\gamma}$.}

In this Section we associate with a  GKZ system another system of  PDE which we call an "antysymmetrized" GKZ system. We define functions $F_{\gamma}$ that form a base in the space of polynomial solutions of this "antysymmetrized" GKZ system.

\subsubsection{ Vectors   $r_i$.} Introduce vectors   $r_1$, $r_2$, $r_3$, written as rows below

\begin{align}
\begin{pmatrix}
a_{-2}&a_{-1}&a_{1}&a_{2}&a_{-2,-1}&a_{-2,1}&a_{-2,2}&a_{-1,1}&a_{-1,2}&a_{1,2}\\
-1&0&1&0&1&0&0&-1&0&0\\
-1&1&0&0&0&0&1&0&-1&0\\
-1&0&1&0&0&0&1&0&0&-1
\end{pmatrix}
\end{align}

To a pair of vectors  $v_1$, $r_1$ there corresponds a Plucker relation

\begin{equation}
\label{s1}
a_{-2}a_{-1,1}-a_{-1}a_{-2,1}+a_{1}a_{-2,-1}=0,
\end{equation}

to a pair of vectors  $v_2$, $r_2$  there corresponds a Plucker relation

\begin{equation}
\label{s2}
a_{2}a_{-2,-1}-a_{-2}a_{2,-1}+a_{-1}a_{-2,2}=0,
\end{equation}

to a pair of vectors   $v_3$, $r_3$ there corresponds a Plucker relation

\begin{equation}
\label{s3}
a_{2}a_{-2,1}-a_{-2}a_{2,1}+a_{1}a_{2,-2}=0,
\end{equation}

One has also a relation

\begin{equation}
\label{s4}
a_{-1,1}=-a_{-2,2}
\end{equation}

\begin{lem}
Any relation between  $a_Y$, $Y\subset \{-2,-1,1,2\}$ is a consequence of relations \eqref{s1}-\eqref{s4}.
\end{lem}

Let us associate with these pairs of vectors some differential operators:  the GKZ operators  $\mathcal{O}_i$, $i=1,2,3$, and their  "antisymmetrizations"   $\bar{\mathcal{O}}_i$, $i=1,2,3$:

\begin{align*}
& \mathcal{O}_1=\frac{\partial^2}{\partial  a_{-2} \partial a_{-1,1}}-\frac{\partial^2}{\partial a_{-1}\partial a_{-2,1}},
\,\,\,\, \bar{\mathcal{O}}_1=\mathcal{O}_1+\frac{\partial^2}{\partial a_{1}\partial a_{-2,-1}},\\
& \mathcal{O}_2=\frac{\partial^2}{\partial a_{2}\partial a_{-2,-1}}-\frac{\partial^2}{\partial a_{-2}\partial a_{2,-1}}
,
\,\,\,\, \bar{\mathcal{O}}_2=\mathcal{O}_2+\frac{\partial^2}{\partial a_{-1}\partial a_{-2,2}},\\
& \mathcal{O}_3=\frac{\partial^2}{\partial a_{2}\partial a_{-2,1}}-\frac{\partial^2}{\partial a_{-2}\partial a_{2,1}}
,
\,\,\,\, \bar{\mathcal{O}}_3=\mathcal{O}_3+\frac{\partial^2}{\partial a_{1}\partial a_{2,-2}}
\end{align*}
 To obtain formulas for the action of $F_{-2,1}$ and  $F_{1,-2}$ we need to obtain a formula for a product of a   $\Gamma$-series and a variable modulo Plucker relations.

This formul is proved in Section  \ref{umn}.  To obtian in in Section   \ref{fnk}  we introduce new classes of functions of hypergeometric type. In   Section  \ref{tchn} using these function we obtain a principle that allows us to prove that a relation holfs modulo  Plucker relations.     Using this principle and th Lemma  \ref{ol}, we finally prove that Lemma  \ref{afg}, which gives a  a formula for a product of a   $\Gamma$-series and a variable modulo Plucker relations.

\subsection{ Functions  $\mathcal{F}_{\gamma}^s$  and   $F_{\gamma}$}
\label{fnk}
 Let $\binom{n}{k}=\frac{n!}{k!(n-k)!}$  be a binomial coefficient.
Let   $X$ denote  an index of a determinant   (thus  $X$ is a subset in $\{-2,-1,1,2\}$) .

Instead of determinants  $a_X$, satisfying the Plucker relations, consider independent variables  $A_X$. Let   $s\in\mathbb{Z}_{\geq 0}^3$, $t\in\mathbb{Z}^3$.
Let us use notations:

$$
sr:=s_1r_1+s_2r_2+s_3r_3,\,\,\, tv:=t_1v_1+t_2v_2+t_3v_3.
$$
 Introduce functions

\begin{equation}
\mathcal{F}_{\gamma}^s(A)=\sum_{t\in\mathbb{Z}^3}\frac{\prod_{i=1}^k\binom{t_i+s_i}{s_i}A^{\gamma-sr+tv}}{(\gamma-sr+tv)!}
\end{equation}

Here we use notations

$$
A^{\gamma-sr+tv}:=\prod_X A_X^{(\gamma-sr+tx)_{\text{ a coordinate with index }X }}
$$

 Introduce functions

$$
F_{\gamma}(A)=\sum_{s\in\mathbb{Z}_{\geq 0}^3}\mathcal{F}_{\gamma}^s(A).
$$

One has immedeately

\begin{equation}
\label{pdf}
\frac{\partial}{\partial A_X}\mathcal{F}_{\gamma}^s(A)=\mathcal{F}_{\gamma-e_X}^s(A),\,\,\,\frac{\partial}{\partial A_X}F_{\gamma}(A)=F_{\gamma-e_X}(A).
\end{equation}

\begin{lem}
Functions  $F_{\gamma}(A)$ are solutions of the system
\begin{equation}\label{systema}\bar{\mathcal{O}}_1F=\bar{\mathcal{O}}_2F=\bar{\mathcal{O}}_3F=0\end{equation}
\end{lem}
\proof

Note that one for binomial coefficients   $\binom{t_i+s_i}{s_i}=\frac{1}{s_i!}(t_i+1)...(t_i+s_i) $ one has

$$
\binom{t_i+s_i}{s_i}-\binom{t_i-1+s_i}{s_i}=\binom{t_i+s_i-1}{s_i-1}
$$

Apply to the function $\mathcal{F}^s_{ \gamma }(A)$ the operator   $\mathcal{O}_1$.   According to the ruler  \eqref{pdf} one gets

\begin{align*}
&\mathcal{O}_1\mathcal{F}^s_{\gamma}(A)=\mathcal{F}^s_{\omega-e_{-2}-e_{-1,1}}(A)-\mathcal{F}^s_{\omega-e_{-1}-e_{-2,1}}(A)=
\\&=\sum_{t\in\mathbb{Z}^3}\frac{(\binom{t_1+s_1}{s_1}-\binom{t_1-1+s_1}{s_1})\prod_{i=2}^3\binom{t_i+s_i}{s_i}A^{\gamma-e_{-2}-e_{-1,1}-sr+tv}}{(\gamma-e_{-2}-e_{-1,1}-sr+tv)!}=\\
&=\sum_{t\in\mathbb{Z}^3}\frac{\binom{t_1+s_1-1}{s_1-1}\prod_{i=2}^3\binom{t_i+s_i}{s_i}A^{\gamma-e_{-2}-e_{-1,1}-sr+tv}}{(\gamma-e_{-2}-e_{-1,1}-sr+tv)!}
\\&=\mathcal{F}^{s-e_1}_{\gamma-e_{-2}-e_{-1,1}}(A)\end{align*}

Hence

$$
\mathcal{O}_1\mathcal{F}^s_{\gamma}(A)=-\frac{\partial^2}{\partial a_{1}\partial a_{-2,-1}}\mathcal{F}^{s-e_1}_{\gamma-r_1}(A)
$$

And thus   $\bar{\mathcal{O}}_1F_{\gamma}(A)=0$.
\endproof


 We call  $F_{\gamma}$ an {\it irreducible } solution of the system  \eqref{systema}

For a monomial  $A^{\gamma}$ we call  $\gamma$  a  {\it support } of this monomial. {\it A support of a function}, written as a sum of a power series is set of support of all its summands.  Denote it as $supp F$.



Take a solution  $F$. Let as represent $suppF$ as a union of the sets of type  $\gamma+B$.
For every such a set take in    $F$  all monomials such that their supports belong to this set.  Denote the resulting functions as  $F^{\gamma}$. If this function satisfies  $\mathcal{O}_{i}(F^{\gamma})=0$, $i=1,2,3$,  then the corresponding  support is called the boundary   (or a  {\it boundary point } in
 $supp F$).   The term point is used because this support really becomes a point if one does all consideration   $mod B$.  Actually we do so.

Obviously an irreducible  $F_{\gamma}$ has a unique boundary point  $\gamma+B$.

\begin{lem}
Every polynomial solution of the system   \eqref{systema} is a sum of irreducible solutions
\end{lem}

\proof Take s solution  $F$ and split it into a sum of functions
$F^{\gamma}$ with supports  $\gamma+B$.

Introduce a partial order on the sets   $\gamma+B$. We say that
$$
\gamma+B\preceq\delta+B,
$$

if  $\gamma+sr=\delta\,\,\, mod B$, $s\in\mathbb{Z}^3_{\geq 0}$.

Since we  are considering only polynomial solution there exist summands $F^{\gamma}$ with supports which are maximal acorrding to this order. Let us show that there supports are are boundary points. Indeed 

$$
\bar{\mathcal{O}}_1F^{\gamma}=\mathcal{O}_1F^{\gamma}+\frac{\partial^2}{\partial a_1\partial a_{-2,-1}}F^{\gamma}
$$

If $suppF^{\gamma}=\gamma+B$, то $supp ( \mathcal{O}_1F^{\gamma})=\gamma-v_1^++B$,  where $v_1^+=e_{-2}+e_{-1,1}$,  and $supp(\frac{\partial^2}{\partial A_1\partial A_{-2,-1}}F^{\gamma})=\gamma-e_1-e_{-2,-1}$. Since  $\bar{\mathcal{O}}_1F=0$, that considering the supports we can conclude that  a summand  $\mathcal{O}_1F^{\gamma}$,  is non-zero then it  must be equal to some of the summands of type  $\frac{\partial^2}{\partial A_1\partial A_{-2,-1}}F^{\delta}$ or  $\mathcal{O}_1(F^{\delta})$ taken with an opposite sign.  Actually for  $\mathcal{O}_1F^{\gamma}$  it is not possible to be equal to a opposite  of  a similar summand but  with another $\delta$, hence  $\mathcal{O}_1F^{\gamma}$  is equal to an opposite of  $\frac{\partial^2}{\partial A_1\partial A_{-2,-1}}F^{\delta}$. Then
$supp F^{\delta}-v^+=\gamma-e_1-e_{-2,-1}$.  This means that $supp F^{\delta}=\gamma+v^+-e_1-e_{-2,-1}+B$.  Hence  $supp F^{\delta}\succeq \gamma+B$,  but the support  $\gamma+B$ is boundary,  hence we obtain a contradiction. Thus  $\mathcal{O}_1F^{\gamma}=0$.

Analogously one proves that  $\mathcal{O}_2F^{\gamma}=\mathcal{O}_3F^{\gamma}=0$.

So our solution has boundary points.
The corresponding functions  $F^{\gamma}$  have supports of type  $\gamma+B$, thus one can write

$$
F^{\gamma}=\sum_{t\in\mathbb{Z}^3}c_t\frac{A^{\gamma+tv}}{(\gamma+tv)!}
$$
 for some number  $c_t$.   Since $F^{\gamma}$
are annihilated by   $\mathcal{O}_i$,  then all  $c_t$ are equal.  Thus
$F^{\gamma}$ are  $\Gamma$-series up to multiplication onto a constant.

Now let us describe a {\bf  procedure}.

\begin{enumerate}
	\item For every boundary point   $\gamma+B$  in  $supp F$ take an irreducible solution  $F_{\gamma}$.
	
	\item  Let us subtract them from
$F$  with such a coefficient that the summands in $F$  with the supports   $\gamma+B$ are reduced.  It is possible since  both in  $F_{\gamma}$ and in   $F$ summands in  $\gamma+B$ form a function which is proportional to a  $\Gamma$-series.

\end{enumerate}

Denote the resulting solution as  $G$.  Let us find boundary points in  $suppG$. These boundary points in  $suppG$ are smaller then the boundary points in   $F$  with respect to the order  $\preceq$.
Let us apply the procedure to   $G$ and so on.

Let us show that  after a finite number of steps  we obtain $0$.  To prove it it is enough to show that  the supports of the resulting functions are contained in some finite set.

For a summand  $F^{\gamma}$ in   $F$ with a maximal support  $\gamma+B$
Let us find the set of those non-negative $s^{\gamma}_{i}$, $i=1,2,3$ such that
$
\gamma-s^{\gamma}_1r_1-s^{\gamma}_2r_2-s^{\gamma}_3r_3+b
$ has only positive coordinates for at least some  $b\in B$. This set is finite. Indeed if one subtracts the vectors   $s^{\gamma}_1r_1$ and  $s^{\gamma}_3r_3$ the the coordinate  $e_1$ reduces,  and this cannot be  compensated by adding   $b\in B$, since it's coordinate  $e_1$ is always $0$. also the substraction of   $e_1$ reduces the coordinate  $e_{-1,1}$. If one subtracts  $s^{\gamma}_2r_2$ then the coordinate  $e_{-1}$ reduces.  This can be compensated only  by adding of the vector $v_1$, but then the coordinate  $e_{-1,1}$ reduces. From these consideration we conclude that   we can  subtract from  $\gamma$ the vectors $r_1,r_2,r_3$ only finite number of times such that $mod B$ we can obtain a vector with positive coordinates.

Introduce notations

$$M_{\gamma}=\bigcup\{\gamma-s^{\gamma}_1r_1-s^{\gamma}_2r_2-s^{\gamma}_3r_2+B\},$$

A union is taken over all  $s_i^{\gamma}$ obtained before.

One has  $suppF_{\gamma}\subset M_{\gamma}$, since $F_{\gamma}=\sum_{s\in\mathbb{Z}^3_{\geq 0}}\mathcal{F}_{\gamma}^s$, and also
$
supp\mathcal{F}_{\gamma}^s=\gamma-s^{\gamma}_1r_1-s^{\gamma}_2r_2-s^{\gamma}_3r_2+B
$,  and a function  $\mathcal{F}_{\gamma}^s$  is non-zero if and only if in it's support there is a  vector with positive coordinates.

One ca easily see that  $\delta+B\prec\gamma+B$ то $M_{\delta}\subset M_{\gamma}$.

From the other hand  $supp F\subset \bigcup_{\gamma} M_{\gamma}$,  where the union is taken over all boundary points  $\gamma$.
Indeed let  $\delta\in suppF$, но $\delta\notin \bigcup_{\gamma} M_{\gamma}$.
Consider $F^{\delta}$.   Analogously to the proof of the fact that maximal  points are boundary points   one can conclude the following. If  $\mathcal{O}_iF^{\delta}\neq 0$,  then  $\delta'=\delta+r_i\in supp F$. Also $\delta+B\prec\delta'+B$ and  $\delta'\notin \bigcup_{\gamma} M_{\gamma}$.  Indeed if  $\delta'\in \bigcup_{\gamma} M_{\gamma}$, then the smaller support alsois contained in this set.
Thus we can increase the support not geting to the set $\bigcup_{\gamma} M_{\gamma}$  untill we obtain a support $\delta''\in supp F$, such that $\mathcal{O}_iF^{\delta''}=0$, $i=1,2,3$. This is a boundary support thus  it belongs to $\bigcup_{\gamma} M_{\gamma}$, we have a constradiction.

Thus at every step of the procedure
{\bf procedure }  the support of the resulting function belongs to   $\bigcup_{\gamma} M_{\gamma}$, where a union is taken over all boundary points of the support  $F$.  This set is finite.  Thus at every step the support reduces and thus  after a finite number of steps we get an empty set.  This means that we represent   $F$  as a sum of functions   $F_{\gamma}$.

\endproof

\subsection{ The main difficulty}
\label{tchn}

The main difficulty in deriving the formulas for the action of generators if the fact that the determinants satisfy some relations. Due to  $a_{-1,1}=-a_{-2,2}$  we can just remove $a_{-2,2}$  and say that the determinant satisfy just the  Plucker relations.

\subsubsection{The main principle}
The key fact is the following notation. The basic Plucker relations are in one-to-one correspondence with operators  $\bar{O}_1,\bar{O}_2,\bar{O}_3$.  

This correspondence leads to the following statement.
If  to a function of determinants  $f(a)$ we collate an operator  $f(\frac{\partial}{\partial A})$, then the following statement holds

\begin{align*}
& \lambda_1f_1(a)+...+\lambda_Nf_N(a)=0 \,\,mod\,\,Plucker\Leftrightarrow\\&\Leftrightarrow  \lambda_1f_1(\frac{\partial}{\partial A})+...+\lambda_Nf_N(\frac{\partial}{\partial A})=0\text{ when acting onto the space of solution of }\bar{\mathcal{O}}_iF=0
\end{align*}

Since the functions  $F_{\omega}$ span the solution space we  formulate this principle asa follows

\begin{lem}
	\label{l4}
	\begin{align*}
	& \lambda_1f_1(a)+...+\lambda_Nf_N(a)=0 \,\,mod\,\,Plucker\Leftrightarrow\\&\Leftrightarrow   \forall \omega\,\,\,(\lambda_1f_1(\frac{\partial}{\partial A})+...+\lambda_Nf_N(\frac{\partial}{\partial A})) F_{\omega}=0
	\end{align*}
\end{lem}

\subsubsection{The main Lemma}
Our main instrument that allows us to obtain a formula for a produc of a $\Gamma$-series and a variable is a formula for the action of   $\mathcal{F}_{\gamma}(\frac{\partial}{\partial A})$ onto  $F_{\omega}(A)$.
\begin{lem}
	\label{ol}
	\begin{equation}
	\label{osnf}
	\mathcal{F}_{\gamma}(\frac{d}{dA})F_{\omega}(A)=\sum_{s\in\mathbb{Z}_{\geq 0}^k}\mathcal{F}_{\gamma+sr}^s(1)F_{\omega-\gamma-sr}(A),
	\end{equation}
	
	where  $\mathcal{F}_{\gamma+sr}^s(1)$ is a result of substitution of  $1$ instead of all arguments
	
\end{lem}

\proof

First of all we need to prove the following relations for the binomial coefficients
\begin{prop}
	$$\binom{N}{a+b}=\sum_{N=N_1+N_2}\binom{N_1}{a}\binom{N_2-1}{b-1}  $$
\end{prop}

\proof
Consider a triangle

$$
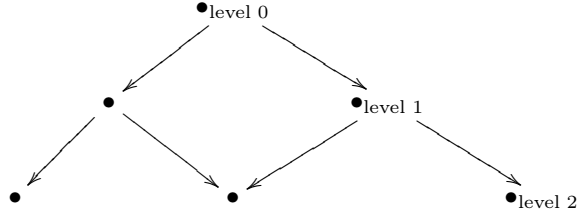
\xymatrix{
	&   & \bullet_{\text{level 0}} \ar[dr]  \ar[dl]             \\
	&\bullet\ar[dl]\ar[dr] & &     \bullet_{\text{level 1}}   \ar[dl]  \ar[dr] \\ \bullet &&\bullet&&\bullet _{\text{level 2}}}
$$

which continues to lower levels.  Then $\binom{N}{a+b}$  is a number of paths from the upper vertex (at level  $0$)  to a vertex at the level $N$,  which has an indent   $a+b$ from the left.  Note on when this path at some level comes to a vertex with an indent   $a$, then until the level  $N_1$ has the same indent,  and then at the level  $N_1+1$  the indent increases. The number of paths which at the level $N_1$ have an indent  $a$, equals to  $\binom{N_1}{a}$.  Then this vertex we take as a beginning.  When we go to the next level we move to the right.  thus the remaining part of the original path gives as a path which at the level  $N_2-1$, where   $N_2=N-N_1$ has an indent  $b-1$. The number of such paths equals to  $\binom{N_2-1}{b-1}$.  To obtain the number of all path satisfying the conditions of Lemma we need to sum over  $N_1$.  Thus we prove the Lemma.
\endproof

\begin{cor}
	\begin{equation}\label{trt}\binom{(t_i+l_i)+k_i}{t_i+l_i}=\sum_{s_i\mathbb{Z}_{\geq 0}}\binom{l_i+s_1-1}{l_i-1}\binom{t_i+k_i-s_i}{t_i-s_i}+...\end{equation}
\end{cor}

Now we return to the proof of  Lemma \ref{ol}.
Let us write  $\mathcal{F}_{\gamma}(\frac{d}{dA})=\sum_l\frac{(\frac{d}{d A})^{\gamma+lv}}{(\gamma+lv)!}$. Let us find an action of  $(\frac{d}{d A})^{\gamma+lv}$ onto a summand  $\mathcal{F}_{\omega}^p(A)$ from   $F_{\omega}$. According to   \eqref{pdf} one has

$$
(\frac{d}{d A})^{\gamma+lv}\mathcal{F}_{\omega}^p(A)=\mathcal{F}_{\omega-\gamma-lv}^p(A).
$$

Consider $\mathcal{F}_{\omega-\gamma-lv}^s(A)$.  We use a notation

$$
\binom{\tau+p}{p}:=\prod_{i=1}^3\binom{\tau_i+p_i}{p_i}.
$$
One has

\begin{align*}
&\mathcal{F}_{\omega-\gamma-lv}^p(A)=\sum_{\tau\in\mathbb{Z}^3}\frac{\binom{\tau+p}{p}A^{\omega-\gamma-pr-lv+\tau v}}{(\omega-\gamma-pr-lv+\tau v)!}=\\
&=\sum_{t\in\mathbb{Z}^3}\frac{\binom{t+l+p}{p}A^{\omega-\gamma-pr+t v}}{(\omega-\gamma-pr+tv)!}
\end{align*}

Apply  \eqref{trt}.  Using  $\sum_{t\in\mathbb{Z}^3}\frac{\binom{t+p-s}{p-s}A^{\omega-\gamma-sr+t v}}{(\omega-\gamma-pr+tv)!}=\mathcal{F}^{p-s}_{\gamma-sr}(A)$,  one gets

$$
\mathcal{F}_{\omega-\gamma-lv}^s(A)=\sum_{s\in\mathbb{Z}^3_{\geq 0}}\binom{l+s-1}{s-1}\mathcal{F}^{p-s}_{\gamma-sr}(A),
$$

where

$$
\binom{l+s-1}{s-1}:=\prod_{i=1}^3\binom{l_i+s_i-1}{s_i-1}
$$

Take an expression for  $(\frac{d}{d A})^{\gamma+lv}\mathcal{F}_{\omega}^p(A)$  and sum them over  $p$,  one gets

$$
(\frac{d}{d A})^{\gamma+lv}F_{\omega}(A)=\sum_{s\in\mathbb{Z}^3_{\geq 0}}\binom{l-1+s}{l-1}F_{\omega-\delta-sr}(A).
$$

Now let us sum over   $l$, one gets

\begin{align*}&\mathcal{F}_{\gamma}(\frac{d}{dA})F_{\omega}(A)=\sum_{s\in\mathbb{Z}^3_{\geq 0}}(\sum_l\frac{\binom{l-1+s}{l-1}}{(\gamma+lv)!})F_{\omega-\delta-sr}(A)=\\
&=\sum_{s\in\mathbb{Z}^3_{\geq 0}}\mathcal{F}^s_{\gamma+v+sr}(1)F_{\omega-\delta-sr}(A)
\end{align*}

\endproof


\subsection{ A formula for a product  of a $A$-hypergeometric function and a variable}
\label{umn}



\begin{lem}
	\label{afg}
	\begin{equation}
	\label{af}
	A_X\mathcal{F}_{\gamma-e_Y}(A)=\sum_p c_p\mathcal{F}_{\gamma-e_Y+e_X+pr}(A) \,\,\,mod\,\,\,Plucker  \,\,\,\,\,
	\end{equation}
	
	where
	
	\begin{align}
	\begin{split}
	\label{cs}
	&c_s=\frac{\mathcal{F}_{\gamma+v-e_Y}^{s}(1)}{\mathcal{F}^s_{\gamma+v-e_Y+e_X+sr}(1)}-\sum_{p=0}^{s-1}\frac{\mathcal{F}_{\gamma+v-e_Y}^{p}(1)\mathcal{F}_{\gamma+v+pr-e_Y+e_X+(s-p)r}^{s-p}(1)}{\mathcal{F}_{\gamma+v+pr-e_Y+e_X+(s-p)r}(1)\mathcal{F}_{\gamma+v+pr-e_Y+e_X}(1)}=\\
	&=
	\frac{\mathcal{F}_{\gamma+v-e_Y}^{s}(1)}{\mathcal{F}^s_{\gamma+v-e_Y+e_X+sr} (1) }-\sum_{p=0}^{s-1}\frac{\mathcal{F}_{\gamma+v-e_Y}^{p}(1)\mathcal{F}_{\gamma+v-e_Y+e_X+sr}^{s-p}(1)}{\mathcal{F}_{\gamma+v-e_Y+e_X+sr}(1)\mathcal{F}_{\gamma+v+pr-e_Y+e_X}(1)}
	\end{split}
	\end{align}
\end{lem}

\proof

To prove the Lemma let us use the principle formulated in  Lemma  \ref{l4},  and the formula  \eqref{osnf}.   Let us transform   \eqref{af} into a differential operator and let us act by this operator onto $F_{\omega}$.  One has

\begin{align*}
&\Big( \frac{d}{dA_X}\mathcal{F}_{\gamma-e_Y}(\frac{d}{dA})\Big )F_{\omega}(A)=\sum_s\mathcal{F}_{\gamma+v-e_Y+sr}^s(1)F_{\omega-\gamma+e_Y-e_X-sr}(A).
\end{align*}

Note that

\begin{align*}
&\mathcal{F}_{\gamma-e_Y+e_X+pr}(\frac{d}{dA})F_{\omega}(A)=\sum_s\mathcal{F}_{\gamma+v-e_Y+e_X+(p+s)r}^s(1)F_{\omega-\gamma+e_Y-e_X-(s+p)r}(A).
\end{align*}

From this formulas one sees that  $A_X\mathcal{F}_{\gamma-e_Y}(A)$ can be expressed though the functions of type  $\mathcal{F}_{\gamma-e_Y+e_X+sr}(A)$  modulo Plucker relations. If

$$
A_X\mathcal{F}_{\gamma-e_Y}(A)=\sum_p c_p\mathcal{F}_{\gamma-e_Y+e_X+pr}(A),
$$

then the coefficients   $c_p$  are solution of the following linear system.  Rows and columns of this system are labed by  $s\in\mathbb{Z}^3_{\geq 0}$. We suppose that  $s\preceq p$, if for all $i$ one has $s_i\leq p_i$. This system is lower-triangular relatively this order

\begin{tiny}
	\begin{equation}
	\label{mtc}
	\begin{pmatrix}
	\mathcal{F}_{\gamma+v-e_Y+e_X}(1) &...& 0 &...&0 &...\\
	...\\
	\mathcal{F}^{s}_{\gamma+v-e_Y+e_X+sr}(1)&...&\mathcal{F}^{}_{\gamma+v-e_Y+e_X+sr}(1) &...&0&...\\
	...\\
	\mathcal{F}^{s+p}_{\gamma+v-e_Y+e_X+(s+p)r}(1)&...&\mathcal{F}^{p}_{\gamma+v-e_Y+e_X+(s+p)r}(1) &...&\mathcal{F}_{\gamma+v-e_Y+e_X+(s+p)r}(1)&...
	
	\end{pmatrix}
	\cdot
	\begin{pmatrix}
	c_0\\...\\c_s\\...
	\end{pmatrix}=
	\begin{pmatrix}
	\mathcal{F}_{\gamma+v-e_Y}(1)\\...\\\mathcal{F}_{\gamma+v-e_Y}^s(1)\\...
	\end{pmatrix}
	\end{equation}
\end{tiny}

One can find explicitly an inverse of the matrix of this system. Note that the firs column of the inverse matrix equals to

$$
(\frac{1}{\mathcal{F}_{\gamma+v-e_Y+e_X}(1)},...,-\frac{\mathcal{F}_{\gamma+v-e_Y+e_X+sr}^s(1)}{\mathcal{F}_{\gamma+v-e_Y+e_X-sr}(1)\mathcal{F}_{\gamma+v-e_Y+e_X}(1)}...)^t
$$

Since instead of  $\gamma$ one can take $\gamma+sr$, then we obtain that the column  $s$ of the matrix of the system  \eqref{mtc} has an analogous form but it begins with  the row   $s$ and we take  $\gamma+sr$   instead $\gamma$.

We conclude that

\begin{align*}
&c_s=\frac{\mathcal{F}_{\gamma+v-e_Y}^{s}(1)}{\mathcal{F}^s_{\gamma+v-e_Y+e_X+sr}(1)}-\sum_{p=0}^{s-1}\frac{\mathcal{F}_{\gamma+v-e_Y}^{p}(1)\mathcal{F}_{\gamma+v+pr-e_Y+e_X+(s-p)r}^{s-p}(1)}{\mathcal{F}_{\gamma+v+pr-e_Y+e_X+(s-p)r}(1)\mathcal{F}_{\gamma+v+pr-e_Y+e_X}(1)}=\\
&=
\frac{\mathcal{F}_{\gamma+v-e_Y}^{s}(1)}{\mathcal{F}^s_{\gamma+v-e_Y+e_X+sr} (1) }-\sum_{p=0}^{s-1}\frac{\mathcal{F}_{\gamma+v-e_Y}^{p}(1)\mathcal{F}_{\gamma+v-e_Y+e_X+sr}^{s-p}(1)}{\mathcal{F}_{\gamma+v-e_Y+e_X+sr}(1)\mathcal{F}_{\gamma+v+pr-e_Y+e_X}(1)}
\end{align*}

\endproof

\subsection{Transformations of diagrams}\label{trs}

To obtain formulas for the action of $F_{-2,1}$ and $F_{1,-2}$ we need  one more calculation.

We   identify a diagram and a shift  vector of a shifted lattice that defiens a    $\Gamma$-series.
Ne need to describe a transformation of a diagram that happens when we add to  $\gamma$ some vectors.

To do it consider the following picture
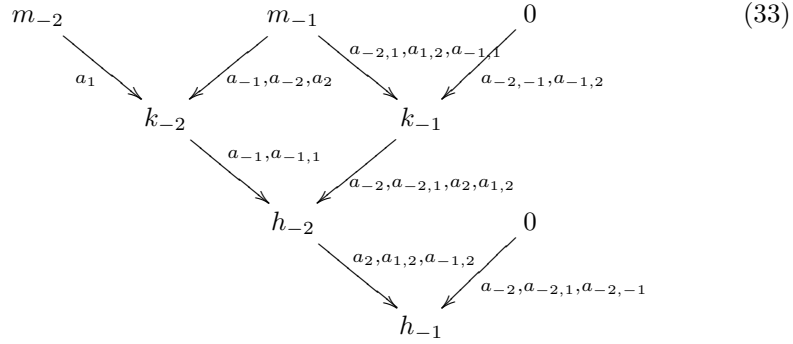
\begin{equation}
\label{dgu}
\xymatrix{
	m_{-2}\ar[dr]_{a_{1}}
	&  &    m_{-1} \ar[dl]^{a_{-1},a_{-2},a_{2}} \ar[dr]^{a_{-2,1},a_{1,2},a_{-1,1}} & & 0  \ar[dl]^{a_{-2,-1},a_{-1,2}}
	\\
	& k_{-2} \ar[dr]^{a_{-1},a_{-1,1}} & & k_{-1}  \ar[dl]^{a_{-2},a_{-2,1},a_{2},a_{1,2}}   \\
	& & h_{-2} \ar[dr]^{a_{2},a_{1,2},a_{-1,2}} &  &  0 \ar[dl]^{a_{-2},a_{-2,1},a_{-2,-1}}\\
	&&& h_{-1}          }
\end{equation}

Equations  \eqref{ur},  that describe the shifted lattice  $\gamma+B$ can be described as follows.  If over an arrow we write some determinantsm that the sum of their powers equals to a difference between numbers occurring at the starting and and the ending of an edge.


Now let us write explicitly transformations of diagrams that happen when we add to   $\gamma$ some vectors.   To obtain formulas  for the action of  $F_{-2,1}$  we add the following vector

$
\begin{array}{|c|c|c|c|c|}
\hline
\text{ a vector which is added to $\gamma$ }&  -e_{1}+e_{-2} & -e_{2}+e_{-1} & -e_{-1,1}+e_{-2,-1} & -e_{1,2}+e_{-1,1} \\
\hline
\text{ a transformation of a diagram }&  \begin{cases}k_{-2}+1,\\ h_{-2}+1,\\ h_{-1}+1\end{cases} & h_{-2}-1  &  \begin{cases} k_{-1}+1 \\h_{-2} +1 \\ h_{-1}+1  \end{cases} &   h_{-2}-1 \\
\hline
\end{array}
$

To obtain formulas  for the action of $F_{1,-2}$ we add the following vector:

$
\begin{array}{|c|c|c|c|c|}
\hline
\text{  a vector which is added  to $\gamma$ }&    -e_{-2}+e_{1} & -e_{-1}+e_{2} & -e_{-2,-1}+e_{-1,1} &-e_{1,-1}+e_{1,2} \\
\hline
\text{ a transformation of a diagram }&  \begin{cases}k_{-2}-1,\\ h_{-2}-1,\\ h_{-1}-1\end{cases} & h_{-2}+1 & \begin{cases} k_{-1}-1 \\h_{-2} -1 \\h_{-1}-1 \end{cases} & h_{-2}+1  \\
\hline
\end{array}
$

Also we investigate the adding of  $r$

$
\begin{array}{|c|c|c|c|}
\hline
\text{ a vector which is added to $\gamma$ }&  r_1 & r_2 & r_3\\
\hline
\text{a transformation of a diagram }&\begin{cases}k_{-2}-1,\\k_{-1}+  1\end{cases} & \begin{cases}k_{-1}-1  \\ h_{-2}-2 \\ h_{-1}-1  \end{cases} &  \begin{cases} k_{-2}-1 \\h_{-2}-2 \\h_{-1}-1 \end{cases}\\
\hline
\end{array}
$

\subsection{ Operators $F_{-2,1}$ and   $F_{1,-2}$.  Formulas for the action}

Apply Lemma  \ref{afg} to the formulas  \eqref{ffd}.  As a result every summand in  \eqref{ffd} is represented as a sum of  $\Gamma$-series with coefficients of type  \eqref{cs}.   A shift vector of these $\Gamma$-series looks as follows   $\gamma-e_Y+e_X+pr$.
The   $\Gamma$-series are the Gelfand-Tselin base vectors, and a transformation of a diagramm corresponding to an adding of   $-e_Y+e_X+pr$ is described in Section \ref{trs}.  Thus we have proved theorems.

\subsubsection{ The action of  $F_{-2,1}$}
\begin{thm}
	The result of an application of  $F_{-2,1}$ is a sum of three series of diagrams. Each series is numbered by    $s^1,s^2,s^2\in\mathbb{Z}_{\geq 0}$.
	
	\begin{enumerate}
		\item The	sum of diagrams obtained from $\gamma$ by transformation
		$$
		\begin{cases}
		k_{-2}-s^1-s^3+1,\,\,\,\,\, k_{-1}+s^1-s^2\\
		h_{-2}-2s^2+2s^3+1\\
		h_{-1}-s^2-s^3+1
		\end{cases}
		$$
		
	each diagram is taken with a coefficient \eqref{cs}, where $e_X=e_{-2}$, $e_Y=e_1$
		\item
		
 The	sum of diagrams obtained from $\gamma$ by transformation
		$$
		\begin{cases}
		k_{-2}-s^1-s^3,\,\,\,\,\, k_{-1}+s^1-s^2\\
		h_{-2}-2s^2+2s^3-1\\
		h_{-1}-s^2-s^3
		\end{cases}
		$$
		
			each diagram is taken with a coefficient , which is obtained as follows. We subtract from the coefficient  \eqref{cs}, where  $e_X=e_{-1}$, $e_Y=e_{2}$  the  double of the coefficient
		\eqref{cs}, where  $e_X=e_{-1,1}$, $e_Y=e_{1,2}$
		\item
		
		The	sum of diagrams obtained from $\gamma$ by transformation
		$$
		\begin{cases}
		k_{-2}-s^1-s^3,\,\,\,\,\, k_{-1}+s^1-s^2+1\\
		h_{-2}-2s^2+2s^3+1\\
		h_{-1}-s^2-s^3+1
		\end{cases}
		$$
		
		each diagram is taken with a minus coefficient \eqref{cs}, where  $e_X=e_{-2,-1}$, $e_Y=e_{-1,1}$

	\end{enumerate}
\end{thm}

\subsubsection{The action of $F_{1,-2}$}
\begin{thm}
		The result of an application of $F_{1,-2}$ is a sum of three series of diagrams. Each series is numbered by    $s^1,s^2,s^2\in\mathbb{Z}_{\geq 0}$.
	
	\begin{enumerate}
		\item 	The	sum of diagrams obtained from $\gamma$ by transformation
		$$
		\begin{cases}
		k_{-2}-s^1-s^3-1,\,\,\,\,\, k_{-1}+s^1-s^2\\
		h_{-2}-2s^2+2s^3-1\\
		h_{-1}-s^2-s^3-1
		\end{cases}
		$$
		
		each diagram is taken with a coefficient \eqref{cs}, where  $e_X=e_{1}$, $e_Y=e_{-2}$
		\item
		
		The	sum of diagrams obtained from $\gamma$ by transformation
		$$
		\begin{cases}
		k_{-2}-s^1-s^3,\,\,\,\,\, k_{-1}+s^1-s^2\\
		h_{-2}-2s^2+2s^3+1\\
		h_{-1}-s^2-s^3
		\end{cases}
		$$
		
	each diagram is taken with a coefficient  \eqref{cs},  where  $e_X=e_{2}$, $e_Y=e_{-1}$,  minus coefficient  \eqref{cs}, where  $e_X=e_{1,2}$, $e_Y=e_{-1,1}$.
		
		\item
		
			The	sum of diagrams obtained from $\gamma$ by transformation
		$$
		\begin{cases}
		k_{-2}-s^1-s^3,\,\,\,\,\, k_{-1}+s^1-s^2-1\\
		h_{-2}-2s^2+2s^3-1\\
		h_{-1}-s^2-s^3-1
		\end{cases}
		$$
		
			each diagram is taken with a minus double of the  coefficient \eqref{cs},  where $e_X=e_{-1,1}$, $e_Y=e_{-2,-1}$

	\end{enumerate}
\end{thm}


\begin{thebibliography}{99}



\bibitem{zh2} D. P. Zhelobenko,  Compact Lie groups and their representations. – American Mathematical Soc., 1973. – V. 40.

\bibitem{M} A.I. Molev, Yangians and Classical Lie Algebras, AMS,  Mathematical Surveys and Monographs, vol. 143, 2007

\bibitem{bb}  G.E.  Biedenharn,  L.C. Baid, On the representations of semisimple Lie Groups II, J. Math. Phys., V. 4, N 12,
1963, 1449-1466.




\bibitem{a1} D.V. Artamonov, Formula for the Product of Gauss Hypergeometric Functions and Applications ,J Math Sci
, 249, 817--826, 2020.

\bibitem{A1}  D.V. Artamonoiv., Clebsh-Gordon coefficients for the algebra  $\mathfrak{gl}_3$ and hypergeometric functions, St. Petersburg Mathematical Journal, accepted

\bibitem{tmph} P. A. Valinevich, “Construction of the Gelfand–Tsetlin basis for unitary principal series representations of the algebra sln(C)”, Theoret. and Math. Phys., 198:1 (2019), 145–155

\bibitem{t15}  V.K. Dobrev, P. Truinin, Polynomial realization of $U_q(sl(3))$ Gel'fand-(Weyl)-Zetlin basis, J. Math. Phys., 38:7 (1997), 3750-3767.

\bibitem{t16}  V.K. Dobrev, A.D.  Mitov, P. Truinin, Normalized $U_q(sl(3))$ Gel'fand-(Weyl)-Zetlin basis and new summation formulas for  $q$-hypergeometric functions, J. Math. Phys., 41:11, 2000, 7752-7768












\bibitem{ap1} N. Hambli, J. Michelson, and R. T. Sharp,  Character states and generator matrix elements for  $Sp(4)\supset U(2)\times U(1)$, Journal of Mathematical Physics 37, 3022 (1996)


\bibitem{ap2} S Alisauskas, Biorthogonal systems for  $SU_4\supset SU_2\times SU_2$, $SU_n\supset SO_n$ and $Sp_4\supset U_2$ and analytical inversion symmetry, J. Phys. A: Math. Gen. , 1987, 20,1045




\bibitem{ap3}  J.A.Evans, N.Kraus, An exact solution of the pairing plus monopole hamiltonian using a boson representation of the group $Sp_4$, Physics Letters B
Volume 37, Issue 5, 27 December 1971, Pages 455-459


\bibitem{ap4} Jin-Quan Chen, Jialun Ping, Fan Wang,  Group representation theory for physicists, World Scientific, 2002

\bibitem{GG} I. M. Gel'fand, M. I. Graev, V. S. Retakh, “General hypergeometric systems of equations and series of hypergeometric type”, Russian Math. Surveys, 47:4 (1992), 1–88



\end{thebibliography}
\end{document}